\documentclass[11pt]{article}
\usepackage{amssymb}

\newtheorem{theorem}{Theorem}

\newtheorem{case}[theorem]{Case}

\input{tcilatex}

\begin{document}

\title{Three Colorability of an Arrangement Graph of Great Circles }
\author{I. Cahit\thanks{%
icahit@ebim.com.tr} \\
Girne American University\\
Girne, North Cyprus}
\date{April 2004}
\maketitle

\begin{abstract}
Stan Wagon asked the following in 2000. Is every zonohedron face 3-colorable
when viewed as a planar map? An equivalent question, under a different
guise, is the following: is the arrangement graph of great circles on the
sphere always vertex 3-colorable?(The arrangement graph has a vertex for
each intersection point, and an edge for each arc directly connecting two
intersection points.) Assume that no three circles meet at a point, so that
this arrangement graph is 4-regular.

In this note we have shown that all arrangement graphs defined as above are
3-colorable.
\end{abstract}

\section{Introduction}

One of the recent interesting graph coloring problem is the following: Is
every zonohedron face 3-colorable when viewed as a planar map [Wag02]? An
equivalent question, under a different guise, is the following: is the
arrangement graph of great circles on the sphere always vertex 3-colorable
[FHNS00] ? Here the arrangement graph has a vertex for each intersection
point, and an edge for each arc directly connecting two intersection points.
Assume that no three circles meet at a point, so that this arrangement graph
is 4-regular. \ Koester has shown that arrangement graphs of circles in the
plane or general circles on the sphere, can require four colors [Koe90]. All
arrangement graphs of up to 11 great circles have been verified to be
3-colorable by O. Aichholzer.

In this note we have shown that all arrangement graphs defined as above are
3-colorable. Our solution is based on the observation that all arrangement
graphs of great circles contain two edge-disjoint identical (actually mirror
image of each other) triangular (cycle of length three) closed chains.

\section{The Triangular Chains}

We assume that the sphere $S$ and the great circles $C_{1},C_{2},...,C_{m}$
are arranged so that one of the great circles say $C_{1}$ is also boundary
at the horizon of $S$ \ when it is seen as bird-eye view. For the sake of
clarity in the figures we have color surface of all triangles on the back
ground of the arrangement graph by dark-grey and color of all triangles on
the fore ground of the arrangement graph on $S$ \ by light-grey. We have the
following simple observations:

(1) Let $G$ be denote the graph corresponding to an arrangement of $k$ great
circles $C_{1},C_{2},...,C_{k}.$ Then the number of triangles created by the
arrangement of $k$ great circles is $2k$.

(2) Let $T_{i}\in S_{1}$ be a fore-ground triangle of $G$ then there exits
mirror image back ground triangle $T_{i}^{\text{,}}\in S_{1}^{,^{\prime }}$
in $G$. For example for $k=4,$ in Fig. 1 mirror image back-ground triangles
of the fore-ground triangles $(ABC),(CEF),(AGH),(BID)$ respectively are

$(A^{\prime }B^{\prime }C^{\prime }),(C^{\prime }HI),(A^{\prime
}DE),(B^{\prime }FD)$.

(3) In the arrangement of great circles graph $G$, the set of $k\geq 4$
triangles can be grouped into two sets as $S_{1}=\{T_{1},T_{2},...,T_{k}\}$
and $S_{1}^{,^{\prime }}=\{T_{1}^{\text{,}},T_{2}^{,},...,T_{k}^{,}\}$ such
that each group of triangles forms a closed triangular edge-disjoint chains
, where $S_{1}\cap S_{1}^{,^{\prime }}=\varnothing $ and $k$ denoting the
number of the great circles. \ 

Furthermore we have $\ \ V(T_{i})\cap V(T_{i+1})=v_{i}$ and $%
V(T_{i}^{^{\prime }})\cap V(T_{i+1}^{\text{'}})=v_{i}^{\prime },i=1,2,...,k$ 
$(\func{mod}k)$ for $k\geq 5.$ For example for $k=4$ in Fig. 1 these
disjoint triangular chains respectively are $\ S_{1}=\left\{
(ABC),(CEF),(FB^{\prime }G),(GHA)\right\} $ and $S_{1}^{\prime }=\left\{
(A^{\prime }B^{\prime }C^{\prime }),(C^{\prime }HI),(IBD),(DEA^{\prime
})\right\} $.

(4) \ Since in any arrangement \ of great circles graph we do not allow more
than two circles to meet at point, for fixed $k$ all great circles graphs
are isomorphic. This will permits us to consider one great circle graph for
the given $k$.

\section{Three-Coloring}

Let us denote three colors by the numbers $1,2,$ and $3$. Our 3-coloring of $%
G$ is based on the disjointedness of the two triangular chains. Depending on
the parity of the number of great circles we have two cases to consider:

\begin{case}
The number of great circles is even.
\end{case}

Again consider the arrangement of four great circles shown in Fig.1. \
Consider the following two cycles:$c_{1}=\left\{ A^{\prime },C^{\prime
},I,D\right\} $ and $c_{2}=\left\{ A,C,F,G\right\} .$ Since these cycles are
disjoint and of even length we can color their vertices with $(1,3)-$(Kempe)
chains. The remaining vertices in $G$ are all the other non-adjacent
vertices of the two triangular chains and therefore three coloring can be
extended \ by coloring them with the color $2,$ i.e., the vertices $%
B,E,B^{\prime },$ and $H$. When the number of great circles is even the
above three coloring of arrangement of four circles can easily be
generalized to any \thinspace $k$ even. For example in Fig. 2 we have shown
three coloring of the arrangement of six great circles. Again the vertices
other than the apex points of the triangles of the two disjoint triangular
closed chains are colored by $(1,3)$-chains. Then apex points of the
triangles in the triangular chains form the \textit{first}-level vertices of 
$G$ and will be colored by the color $2$. Similarly vertices (not yet
colored) adjacent to the first level vertices that have been colored by $2$
form the \textit{second-}level vertices of $G$ and we color them by the
color $1$. Actually the first and second level vertices can be viewed as the
points of a rectangular-lattice $L_{r\text{ }}$points embedded on the
surface of the sphere $S$. For $k=6$ the number of the lattice levels is two
and it is enough to color these level-vertices respectively by colors $1$
and $2$ (see Fig. 2). In general for $k$ great circles arrangement graph $G$
the number of the rectangular lattice levels will be $k-4\,\ $\ and the
vertices on the crossing points of the lattice are colored alternatingly
with the colors $\ 1$ and $2$. Note that, since the two triangular closed
chains that have been colored by $(1,3)$-chains, the overall coloring of $G$
is a proper $3$ coloring.

\begin{case}
The number of great circles is odd.
\end{case}

This case is different than the above since the length of the two disjoint
triangular chains are odd. We will demonstrate our three coloring method for
the arrangement graphs with $k=5,7$ (see Fig. 3,4) and explain the details
for any odd $k$.

$k\equiv 0(\func{mod}2)$ : Consider the sequence of the colors in the two
triangular chains in Fig.3, that is $T_{1}=\left\{ 3^{\ast
},1,2,1,3,2,3^{\ast }\right\} $ and $\ T_{2}=\left\{
1^{\#},3,2,3,1,2,1^{\#}\right\} ,$ where "*" and "$\#$" denote the color of
the apex vertices of the light and dark-grey triangles. The coloring of the
other vertices which are the apex points of the triangles is shown in Fig. 3.

$k\equiv 1(\func{mod}2):$ Consider the sequence of the colors in the two
triangular chains in Fig. 4, that is $T_{1}=\left\{ 3^{\ast
},2,3,2,3,1,3,1,3^{\ast }\right\} $ and $T_{2}=\left\{
3^{\#},1,3,2,3,1,3,2,3^{\#}\right\} $, where "$\ast $" pointing the color of
the apex-vertex in the light-grey triangle and "$\#$" pointing the color of
the apex vertex in the dark-grey triangle in the triangular chains of $G$. \
Starting these 3-colorings we can extend it to the other vertices of the
graph as follows, where arrows show propagation of coloring on the vertices
of $G$ as:

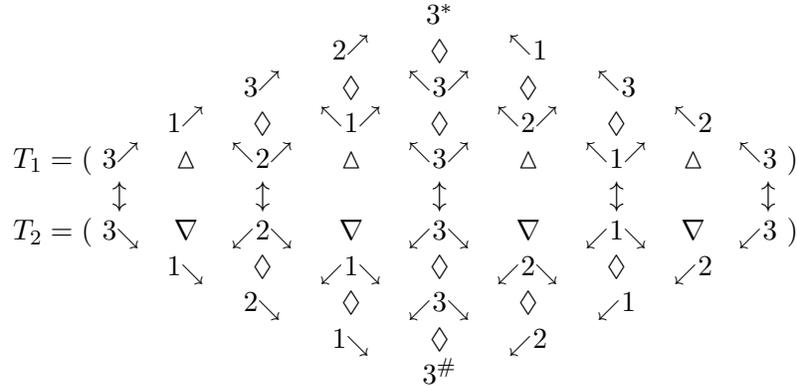
\begin{center}
$%
\begin{array}{ccccccccc}
&  &  &  & 3^{\ast } &  &  &  &  \\ 
&  &  & 2^{\nearrow } & {\huge \Diamond } & ^{\nwarrow }1 &  &  &  \\ 
&  & 3^{\nearrow } & {\huge \Diamond } & ^{\nwarrow }3^{\nearrow } & {\huge %
\Diamond } & ^{\nwarrow }3 &  &  \\ 
& 1^{\nearrow } & {\huge \Diamond } & ^{\nwarrow }1^{\nearrow } & {\huge %
\Diamond } & ^{\nwarrow }2^{\nearrow } & {\huge \Diamond } & ^{\nwarrow }2 & 
\\ 
T_{1}=(\text{ }3^{\nearrow } & {\Huge \vartriangle } & ^{\nwarrow
}2^{\nearrow } & {\Huge \vartriangle } & ^{\nwarrow }3^{\nearrow } & {\Huge %
\vartriangle } & ^{\nwarrow }1^{\nearrow } & {\Huge \vartriangle } & 
^{\nwarrow }3\text{ }) \\ 
\text{ \ \ \ \ \ \ \ }\updownarrow &  & \updownarrow &  & \updownarrow &  & 
\updownarrow &  & \updownarrow \\ 
T_{2}=(\text{ }3_{\searrow } & {\LARGE \nabla } & _{\swarrow }2_{\searrow }
& {\LARGE \nabla } & _{\swarrow }3_{\searrow } & {\LARGE \nabla } & 
_{\swarrow }1_{\searrow } & {\LARGE \nabla } & _{\swarrow }3\text{ }) \\ 
& 1_{\searrow } & {\huge \Diamond } & _{\swarrow }1_{\searrow } & {\huge %
\Diamond } & _{\swarrow }2_{\searrow } & {\huge \Diamond } & _{\swarrow }2 & 
\\ 
&  & 2_{\searrow } & {\huge \Diamond } & _{\swarrow }3_{\searrow } & {\huge %
\Diamond } & _{\swarrow }1 &  &  \\ 
&  &  & 1_{\searrow } & {\huge \Diamond } & _{\swarrow }2 &  &  &  \\ 
&  &  &  & 3^{\#} &  &  &  & 
\end{array}%
$
\end{center}

Symmetric of the coloring of the $T_{1}$ and $T_{2}$ have been shown by "$%
\updownarrow "$ in the Fig. 4. Color the other vertices which are on the
circle $C_{1}$ i.e., great circle at the horizon of sphere $S,$ as $\left\{
1,2,1,3,1,2,1,2,3,2,1,2,1\right\} $.

\bigskip

For the other odd values of $k$ the sequence of colors of the triangular
chains shown above can be generalized easily.

\FRAME{ftbpFU}{289.375pt}{203.3125pt}{0pt}{\Qcb{Four great circles of
arrangement graph.}}{}{Figure}{\special{language "Scientific Word";type
"GRAPHIC";maintain-aspect-ratio TRUE;display "USEDEF";valid_file "T";width
289.375pt;height 203.3125pt;depth 0pt;original-width
419.125pt;original-height 293.9375pt;cropleft "0";croptop "1";cropright
"1";cropbottom "0";tempfilename 'HZYNM401.wmf';tempfile-properties "XPR";}}

\begin{center}
\FRAME{ftbpFU}{284.25pt}{236pt}{0pt}{\Qcb{Three coloring of the six-great
circles graph.}}{}{Figure}{\special{language "Scientific Word";type
"GRAPHIC";maintain-aspect-ratio TRUE;display "USEDEF";valid_file "T";width
284.25pt;height 236pt;depth 0pt;original-width 673.5625pt;original-height
558.625pt;cropleft "0";croptop "1";cropright "1";cropbottom "0";tempfilename
'I09SHK00.wmf';tempfile-properties "XPR";}}\FRAME{ftbpFU}{305.1875pt}{220pt}{%
0pt}{\Qcb{Three coloring of the five-great circles graph.}}{}{Figure}{%
\special{language "Scientific Word";type "GRAPHIC";maintain-aspect-ratio
TRUE;display "USEDEF";valid_file "T";width 305.1875pt;height 220pt;depth
0pt;original-width 461.8125pt;original-height 332.4375pt;cropleft
"0";croptop "1";cropright "1";cropbottom "0";tempfilename
'I0QHU501.wmf';tempfile-properties "XPR";}}\FRAME{ftbpFU}{299.1875pt}{%
236.75pt}{0pt}{\Qcb{Three coloring of the seven great circles graph.}}{}{%
Figure}{\special{language "Scientific Word";type
"GRAPHIC";maintain-aspect-ratio TRUE;display "USEDEF";valid_file "T";width
299.1875pt;height 236.75pt;depth 0pt;original-width
725.3125pt;original-height 573.125pt;cropleft "0";croptop "1";cropright
"1";cropbottom "0";tempfilename 'I0N2UI00.wmf';tempfile-properties "XPR";}}%
\bigskip 
\end{center}

\bigskip

\bigskip

\bigskip

\end{document}